\documentclass{amsproc}

\usepackage[dvips]{graphicx}
\usepackage{amsmath,amssymb}

\begin{document}

\title[Inverse dynamic problem for canonical systems.]
{Inverse dynamic problems for canonical systems  and de Branges
spaces.}

\author[A.\,S.~Mikhaylov, V.\,S.~Mikhaylov]
{$^1$, $^2$A.\,S.~Mikhaylov, $^1$, $^2$V.\,S.~Mikhaylov}

\address{
$^1$ St. Petersburg Department of V.A. Steklov Institute of
Mathematics of the Russian Academy of Sciences, 7, Fontanka,
191023 St. Petersburg, Russia. $^2$ Saint Petersburg State
University, St.Petersburg State University, 7/9 Universitetskaya
nab., St. Petersburg, 199034 Russia.}

\email{mikhaylov@pdmi.ras.ru, vsmikhaylov@pdmi.ras.ru}


\begin{abstract}
We show the equivalence of inverse problems for different
dynamical systems and corresponding canonical systems. For
canonical system with general Hamiltonian we outline the strategy
of studying the dynamic inverse problem and procedure of
construction of corresponding de Branges space.
\end{abstract}


\maketitle

\newtheorem{definition}{Definition}
\newtheorem{lemma}{Lemma}
\newtheorem{proposition}{Proposition}
\newtheorem{theorem}{Theorem}

\section{Introduction}

This is an accompanying paper to \cite{MM2}, in which the authors
have shown the relationship between the de Branges method and the
Boundary Control (BC) method on a basis of three dynamical
systems: wave equation with a potential on a half-line, Dirac
system on a half-line and dynamical system with discrete time for
semi-infinite discrete Schr\"odinger operator. For each system
they constructed the related de Branges space using natural
dynamic objects and operators, used in the BC method. In the
present note we will show the equivalence of dynamic inverse
problems (IP) for different dynamical systems (wave equation,
Dirac system, Jacobi matrices), and IPs for equivalent canonical
systems. We note that every original system will be equivalent to
canonical system with different dynamics (the dependence on $t$ is
given by one of the following operators:
$\frac{d^2}{dt^2},\,i\frac{d}{dt},\partial_t$, where $\partial_t$
is a difference operator).

Let $H\in L_{1,loc}(0,L;\mathbb{R}^{2\times2})$ be a locally
summable on $(0,L)$, $L\leqslant\infty$ matrix-valued function
$H\geqslant 0$,
called Hamiltonian, $J:=\begin{pmatrix} 0 & 1\\
-1 & 0\end{pmatrix}$,  vector $Y=\begin{pmatrix} Y_1 \\
Y_2\end{pmatrix}$. We choose the "proper" dynamics and fix the
general dynamical canonical system, the initial boundary value
problem (IBVP) of which will be the subject of our interest:
\begin{equation*}
iH\frac{dY}{dt}-J\frac{dY}{dx}=0,\quad x\geqslant 0,\, t\geqslant
0.
\end{equation*}
For such a system we set up an IP and outline the strategy of
solving it by the BC method, provided the Hamiltonian is smooth
and strictly positive. We also provide a method of construction of
the de Branges space for such a Hamiltonian in natural dynamic
terms following \cite{MM2}.

In the second section we expose all necessary information on de
Branges spaces and canonical systems following \cite{R2} and
\cite{RR}. In the third section we deal with dynamical systems for
Schr\"odinger operator on a half-line, wave equation on a
half-line, Dirac operator on a half-line and a semi-infinite
Jacobi matrices. We formulate dynamic IP for each system, then we
transform IBVP for each system to the IBVP for certain canonical
system, formulate IP for canonical system, and show that it is
equivalent to original ones.

In the forth section we will show that one specific choice of
dynamics give a finite speed of wave propagation in a canonical
system, provided the Hamiltonian is smooth and strictly positive.
We note that the finiteness of the wave propagation is important:
initially the BC method was developed and applied in the case of
multidimensional wave equation \cite{B87,B97} on a bounded
manifold, but later on the BC method was successfully applied to
parabolic and Schr\"odinger equations (where the speed is
infinite) as well\cite{BIP07,B17,ALP1}. We provide algorithms of
solving dynamic IP and construction of de Branges space for such a
Hamiltonian. Basing on these results, we formulate the hypothesis
on construction the de Branges space for general Hamiltonian by
the dynamic method.

\section{de Branges spaces.}

Here we provide the information on de Branges spaces in accordance
with \cite{R2,RR}. The entire function $E:\mathbb{C}\mapsto
\mathbb{C}$ is called a \emph{Hermite-Biehler function} if
$|E(z)|>|E(\overline z)|$ for $z\in \mathbb{C}_+$. We use the
notation $F^\#(z)=\overline{F(\overline{z})}$. The \emph{Hardy
space} $H_2$ is defined by: $f\in H_2$ if $f$ is holomorphic in
$\mathbb{C}^+$ and
$\sup_{y>0}\int_{-\infty}^\infty|f(x+iy)|^2\,dx<\infty$. Then the
\emph{de Branges space} $B(E)$ consists of entire functions such
that:
\begin{equation*}
B(E):=\left\{F:\mathbb{C}\mapsto \mathbb{C},\,F \text{ entire},
\int_{\mathbb{R}}\left|\frac{F(\lambda)}{E(\lambda)}\right|^2\,d\lambda<\infty,\,\frac{F}{E},\frac{F^\#}{E}\in
H_2\right\}.
\end{equation*}
The space $B(E)$ with the scalar product
\begin{equation*}
[F,G]_{B(E)}=\frac{1}{\pi}\int_{\mathbb{R}}\overline{ F(\lambda)}
G(\lambda)\frac{d\lambda}{|E(\lambda)|^2}
\end{equation*}
is a Hilbert space. For any $z\in \mathbb{C}$ the
\emph{reproducing kernel} is introduced by the relation
\begin{equation}
\label{repr_ker} J_z(\xi):=\frac{\overline{E(z)}E(\xi)-E(\overline
z)\overline{E(\overline \xi)}}{2i(\overline z-\xi)}.
\end{equation}
Then
\begin{equation*}
F(z)=[J_z,F]_{B(E)}=\frac{1}{\pi}\int_{\mathbb{R}}\overline{J_z(\lambda)}
F(\lambda)\frac{d\lambda}{|E(\lambda)|^2}.
\end{equation*}
We observe that a Hermite-Biehler function $E(\lambda)$ defines
$J_z$ by (\ref{repr_ker}). The converse is also true
\cite{DMcK,DBr}: a Hilbert space of analytic functions with
reproducing kernel is a de Branges space (provided some
nonrestrictive conditions on the set of function and on the norm
hold true).


Let $H\in L_{1,loc}(0,L;R^{2\times2})$ be a Hamiltonian and the vector $Y=\begin{pmatrix} Y_1 \\
Y_2\end{pmatrix}$ be solution to the following Cauchy problem:
\begin{eqnarray}
\label{CanSyst}-J\frac{dY}{dx}=\lambda HY,\\
Y(0)=C,\notag
\end{eqnarray}
for $C\in \mathbb{R}^2$, $C\not=0$. Without loss of generality it
is assumed that $\operatorname{tr}H(x)=1.$ Then the function
$E_x(\lambda)=Y_1(x,\lambda)+iY_2(x,\lambda)$ is a Hermite-Biehler
function ($E_L(\lambda)$ makes sense if $L<\infty$), it is called
de Branges function of the system (\ref{CanSyst}) since one can
construct de Branges space based on this function. On the other
hand, $E_L$ serves as an inverse spectral data for the canonical
system (\ref{CanSyst}). The solution to (\ref{CanSyst}) and
$Y(0)=(1,0)^T$ is denoted by $\Theta(x,\lambda)$. The main result
of the theory \cite{RR,DBr} says that the opposite is also true:
every Hermite-Biehler function satisfying some condition comes
from some canonical system.

\section{Dynamical canonical systems for wave equation, Dirac system and Jacobi system with discrete time.}

In this section we use some ideas from \cite{RR} to rewrite IBVPs
for different dynamical systems as IBVPs for canonical dynamical
systems. Everywhere below, $T>0$ is fixed.

\subsection{Wave equation with a potential on a half-line.}

For a potential $q\in L_{1,\operatorname{loc}}(\mathbb{R}_+)$ we
consider the IBVP for the 1d wave equation on a half-line:
\begin{equation}
\label{wave_eqn} \left\{
\begin{array}l
u_{tt}(x,t)-u_{xx}(x,t)+q(x)u(x,t)=0, \quad x\geqslant 0,\ t\geqslant 0,\\
u(x,0)=u_t(x,0)=0,\ u(0,t)=f(t).
\end{array}
\right.
\end{equation}
Here $f$ is an arbitrary $L^2_{loc}\left( \mathbb{R}_{+}\right) $
function referred to as a \emph{boundary control}. The
\emph{response operator} $R^T_q: L_2(0,T)\mapsto L_2(0,T)$ with
the domain $\mathcal{D}=C_0^\infty(0,T)$ is introduced by
$\left(R_q^Tf\right)(t):=u^f_x(0,t)$, it plays a role of a dynamic
inverse data \cite{AM,BM_Sh,MM3}. The IP is to recover q on
$(0,T)$ from $R_q^{2T}$.

We consider the solutions $y_{1,2}$ to following Cauchy problems:
\begin{equation}
\label{Schr_Cauchy} \left\{
\begin{array}l
-y_{1,2}''(x)+q(x)y_{1,2}(x)=0, \quad x\geqslant 0,\\
y_1(0)=1, \, y'_1(0)=0,\, y_2(0)=0, y_2'(0)=1,
\end{array}
\right.
\end{equation}
and look for the solution to (\ref{wave_eqn}) in a form
\begin{equation}
\label{U_repr}
u^f(x,t)=c^1(x,t)y_1(x)+c^2(x,t)y_2(x).
\end{equation}
Plugging this representation to (\ref{wave_eqn}) yields:
\begin{eqnarray*}
c^1_{tt}y_1+c^2_{tt}y_2=-qc^1y_1-qc^2y_2+c^1_{xx}y_1+2c^1_xy_1'+c_1y_1''+c^2_{xx}y_2+2c^2_xy_2'+c_2y''\\
=\left(c^1_xy_1+c^2_xy_2\right)_x+c^1_xy_1'+c^2_xy_2'.
\end{eqnarray*}
If we demand the equality $c^1_xy_1+c^2_xy_2=0,$ then unknown
$c^{1,2}$ satisfy the following system:
\begin{equation}
\label{S1} \left\{\begin{array}l
c^1_{tt}y_1+c^2_{tt}y_2=c^1_xy_1'+c^2_xy_2',\\
c^1_xy_1+c^2_xy_2=0.
\end{array}\right.
\end{equation}
We note that due to the boundary conditions in (\ref{Schr_Cauchy})
and (\ref{S1}), we have that
\begin{equation*}
u^f_x(0,t)=c^1_x(0,t)y_1(0)+c^1(0,t)y_1'(0)+c^2_x(0,t)y_2(0)+c^2(0,t)y_2'(0)=c^2(0,t).
\end{equation*}
On expressing $c^{1,2}_x$ from (\ref{S1}),  and bearing in mind
the equality $\det{\begin{pmatrix}y_1 & y_2\\ y_1' &
y_2'\end{pmatrix}}=1$, we obtain that
\begin{equation*}
\left\{\begin{array}l
c^1_x=-c^1_{tt}y_1y_2-c^2_{tt}y_2^2,\\
c^2_x=c^1_{tt}y_1^2+c^2_{tt}y_1y_2.
\end{array}\right.
\end{equation*}
On introducing the notations $C=\begin{pmatrix}
c^1\\c^2\end{pmatrix}$, $J=\begin{pmatrix} 0& 1\\ -1&
0\end{pmatrix},$ $H=\begin{pmatrix}y_1^2 & y_1y_2\\
y_1y_2& y_2^2\end{pmatrix}$ and counting the initial and boundary
conditions on $u$ at $t=0$ and at $x=0$, we obtain that $C$
satisfies the following IBVP:
\begin{equation}
\label{Can_Syst_Schr} \left\{\begin{array}l
HC_{tt}-JC_x=0,\quad x\geqslant 0,\,t\geqslant 0,\\
C(x,0)=0, C_t(x,0)=0,\quad x\geqslant 0,\\
c^1(0,t)=f(t),\quad t\geqslant 0.
\end{array}\right.
\end{equation}
The response operator $\widetilde R^T_q: L_2(0,T) \mapsto
L_2(0,T)$ for (\ref{Can_Syst_Schr}) is introduced by the equality
$\left(\widetilde R^T_sf\right)(t):=c^2(0,t)$. On the other hand,
using (\ref{U_repr}) and second line in (\ref{S1}), we have that
\begin{equation*}
\left(R_q^Tf\right)(t):=u^f_x(0,t)=c^1_x(0,t)y_1(0)+c^2(0,t)y_2'(0)=c^2(0,t)=\left(\widetilde
R^T_sf\right)(t).
\end{equation*}
So we can see that IPs for (\ref{wave_eqn}) and for
(\ref{Can_Syst_Schr}) are equivalent.

\subsection{Wave equation on a half-line.}

For a smooth positive density $\rho\in C^2(\mathbb{R}_+)$,
$\rho(x)\geqslant\delta>0$, we consider the IBVP for a wave
equation on a half-line:
\begin{equation}
\label{wave_eqn1} \left\{
\begin{array}l
\rho(x)u_{tt}(x,t)-u_{xx}(x,t)=0, \quad x\geqslant 0,\ t\geqslant 0,\\
u(x,0)=u_t(x,0)=0,\ u(0,t)=f(t).
\end{array}
\right.
\end{equation}
Where the function $f\in L^2_{loc}\left(
\mathbb{R}_{+},\mathbb{C}\right) $ is interpreted as a
\emph{boundary control}. The \emph{response operator} $R^T_\rho:
L_2(0,T)\mapsto L_2(0,T)$ with the domain
$\mathcal{D}=C_0^\infty(0,T)$ is defined by $R^T_\rho
f:=u^f_x(0,t)$. We introduce the \emph{eikonal}
$\tau(x):=\int_0^x\rho^{\frac{1}{2}}(s)\,ds$, from physical point
of view, it is a time at which a wave initiated at $x=0$ fills the
segment $(0,x)$, let $\Omega^l=\{x>0\,|\, \tau(x)<l\}$. Then the
natural set up of IP is to recover $\rho(x)|_{\Omega^T}$ from
$R^{2T}_\rho$, see \cite{B08}.

We introduce the new function
\begin{equation*}
C(x,t)=\begin{pmatrix}c^1\\ c^2\end{pmatrix}:=\begin{pmatrix}u_t\\
iu_x\end{pmatrix},
\end{equation*}
and a Hamiltonian $H:=\begin{pmatrix} \rho(x) & 0\\ 0&
1\end{pmatrix}$. Then it is easy to see that $Y$ satisfies the
canonical system
\begin{equation}
\label{Wave_Can_Syst} \left\{
\begin{array}l
iHC_t-JC_x=0,\quad x\geqslant 0,\,t\geqslant 0,\\
C(x,0)=0,\quad x\geqslant 0,\\
c^1(0,t)=g(t):=f'(t), \quad t\geqslant 0.
\end{array}
\right.
\end{equation}
The response operator $\widetilde R^T_\rho: L_2(0,T) \mapsto
L_2(0,T)$ for (\ref{Wave_Can_Syst}) with the domain
$\mathcal{D}=C_0^\infty(0,T)$ is introduced by $\left(\widetilde
R^T_sg\right)(t):=c^2(0,t)$. We can see that IPs for
(\ref{wave_eqn1}) and  for (\ref{Wave_Can_Syst}) are equivalent.

\subsection{Dirac system on a half-line.}

With a matrix potential $V=\begin{pmatrix} p&q\\q&-p
\end{pmatrix}$, $p,q\in C^1_{loc}(\mathbb{R}_{+})$,  vector $u=\begin{pmatrix} u_1 \\
u_2\end{pmatrix}$ we associate the IBVP for a Dirac system:
\begin{equation}
\label{d1} \left\{
\begin{array}l
iu_t+Ju_x+Vu=0,\quad  x\geqslant 0,\,t\geqslant 0,\\
u\big|_{t=0}=0, \quad x\geqslant 0, \\
u_1\big|_{x=0}=f,\,\,  t\geqslant 0,
\end{array}
\right.
\end{equation}
Here $f$ is an arbitrary $L^2_{loc}\left(
\mathbb{R}_{+},\mathbb{C}\right) $ function referred to as a
\emph{boundary control}. The \emph{response operator} $R^T_D:
L_2(0,T)\mapsto L_2(0,T)$ with the domain
$\mathcal{D}=C_0^\infty(0,T)$ is introduced by
$\left(R_D^Tf\right)(t):=u_2(0,t)$, it plays a role of a dynamic
inverse data. The IP is to recover $V$ on $(0,T)$ from $R_D^{2T}$,
see \cite{BM_Dir}.

Let $Y^{1,2}$ be solutions to the following Cauchy problems:
\begin{equation*}
\left\{
\begin{array}l
JY^{1,2}_x+VY^{1,2}=0,\\
Y^1_1(0)=1,\,Y^1_2(0)=0,\, Y^2_1(0)=0,\, Y^2_2(0)=1.
\end{array}
\right.
\end{equation*}
We will look for the solution to (\ref{d1}) in the form
\begin{equation}
\label{U_dir_repr} u(x,t)=c^1(x,t)Y^1(x)+c^2(x,t)Y^2(x)
\end{equation}
Plugging this representation in (\ref{d1}) yields
\begin{eqnarray*}
i\left(c^1_tY^1+c^2_tY^2\right)+c^1_xJY^1+c^2_xJY^2+c_1JY^1_x+c_2JY^2_x+c_1VY^1+c_2VY^2\\
=i\left(c^1_tY^1+c^2_tY^2\right)+J\left(c^1_xY^1+c^2_xY^2\right)=0,
\end{eqnarray*}
on introducing $C=\begin{pmatrix} c^1 \\ c^2\end{pmatrix}$, we see
that the above equality is equivalent to
\begin{equation*}
i\begin{pmatrix}Y^1_1 & Y^2_1 \\ Y^1_2 &
Y^2_2\end{pmatrix}C_t+J\begin{pmatrix}Y^1_1 & Y^2_1 \\ Y^1_2 &
Y^2_2\end{pmatrix} C_x=0.
\end{equation*}
Introduce the notation: $A=\begin{pmatrix}Y^1_1 & Y^2_1 \\ Y^1_2 &
Y^2_2\end{pmatrix}$, $B=JAJ$. Then the above system is equivalent
to
\begin{equation*}
iAC_t-BJC_x=0,
\end{equation*}
on multiplying it by $B^{-1}$ and introducing the Hamiltonian by
$H=B^{-1}A$,  we obtain
\begin{equation*}
iHC_t-JC_x=0,
\end{equation*}
Counting that $\det{B}=\det{A}=1$, we evaluate:
\begin{equation*}
H=B^{-1}A=\begin{pmatrix} Y^1Y^1 & Y^1Y^2\\ Y^1Y^2 &
Y^2Y^2\end{pmatrix},
\end{equation*}
Bearing in mind the initial and boundary conditions in (\ref{d1}),
we see that $C$ satisfies the following IBVP:
\begin{equation}
\label{Dirac_Can_Syst}
\left\{
\begin{array}l
iHC_t-JC_x=0,\quad x\geqslant 0,\,t\geqslant 0,\\
C(x,0)=0,\quad x\geqslant 0,\\
c^1(0,t)=f(t)\quad t\geqslant 0.
\end{array}
\right.
\end{equation}
The response operator $\widetilde R^T_D: L_2(0,T) \mapsto
L_2(0,T)$ for (\ref{Dirac_Can_Syst}) is introduced by
$\left(\widetilde R^T_Df\right)(t):=c^2(0,t)$. The representation
(\ref{U_dir_repr}) implies that IPs for (\ref{d1}) and for
(\ref{Dirac_Can_Syst}) are equivalent.

\subsection{Semi-infinite Jacobi matrices.}

Let $0=b_0<b_1<b_2<\ldots<b_n<\ldots$ be a partition of
$[0,+\infty)$. We introduce the notations:
$\Delta_j:=(b_{j-1},b_j),$ $l_j=|\Delta_j|=b_j-b_{j-1}.$ Let for
each $j$ we define $e_j\in \mathbb{R}^2$, $|e_j|=1$, $e_j\not=\pm
e_{j\pm1},$ and $e_j(x)=e_j$, $x\in\Delta_j.$ We define a
Hamiltonian $H$:
\begin{equation*}
H(x)f(x)=\left(f(x),e_j(x)\right)e_j(x)=\begin{pmatrix}e^2_{1j}(x)
& e_{1j}(x)e_{2j}(x)\\ e_{1j}(x)e_{2j}(x) &
e^2_{2j}(x)\end{pmatrix}\begin{pmatrix}f^1(x) \\ f^2(x)
\end{pmatrix}
\end{equation*}
Consider functions of the type (i.e. functions from the domain of
operator, corresponding to such a Hamiltonian, see [RR]).
\begin{equation}
\label{F_form}
f(x)=\begin{pmatrix}f^1(x)\\f^2(x)\end{pmatrix}=f_je_j(x)+\xi_j(x)e_j^\bot(x),\quad
x\in \Delta_j, f_j\in \mathbb{R},\quad e_j^\bot=Je_j
\end{equation}
and note that $(f,e_j)=f_j.$ For such a Hamiltonian $H$ we study
the equation
\begin{equation}
\label{EqnCS} Jf'=Hg,
\end{equation}
where the function $g$ has a form (\ref{F_form}),
$g=g_je_j(x)+\eta_j(x)e_j^\bot(x)$, $x\in \Delta_j$. The equality in
(\ref{EqnCS}) implies that
\begin{equation*}
\xi'_j(x)Je_j^\bot(x)=g_je_j(x), \quad x\in \Delta_j,
\end{equation*}
which yields the following expression for $\xi_j(x)$ for some $s_j$:
\begin{equation}
\label{Xi_expr}
 \xi_j(x)=s_j+g_j(b_j-x),\quad x\in \Delta_j.
\end{equation}
We use the continuity condition at $x=b_{j-1}$ to get
\begin{equation*}
f_{j-1}e_{j-1}+s_{j-1}e_{j-1}^\bot=f_{j}e_{j}+\left(s_{j}+g_jl_j\right)e_{j}^\bot.
\end{equation*}
Multiplying the above equality by $e_j$ we get
\begin{equation}
\label{f_1}
s_{j-1}=\frac{1}{\left(e_j,e_{j-1}^\bot\right)}\left(f_{j}-f_{j-1}\left(e_j,e_{j-1}^\bot\right)\right),
\end{equation}
and multiplying by $e_{j-1}$ we obtain
\begin{equation}
\label{f_2}
f_{j-1}=f_{j}\left(e_j,e_{j-1}\right)+\left(s_j+g_jl_j\right)\left(e_j^\bot,e_{j-1}\right).
\end{equation}
Using (\ref{f_1}), (\ref{f_2}) we can express $g_j$ via $f_{j-1}$,
$f_{j}$, $f_{j+1}$:
\begin{equation}
\label{f_3}
g_jl_j=\frac{1}{\left(e_j,e_{j-1}^\bot\right)}f_{j-1}+\left(\frac{\left(e_{j+1},e_{j}\right)}{\left(e_{j+1},e_{j}^\bot\right)}-\frac{\left(e_{j},e_{j-1}\right)}{\left(e_{j}^\bot,e_{j-1}\right)}\right)f_j
-\frac{1}{\left(e_j^\bot,e_{j+1}\right)}f_{j+1}.
\end{equation}
Making the substitution
\begin{equation}
\label{Subst} u_j=g_j\sqrt{l_j},\quad v_j=f_j\sqrt{l_j},
\end{equation}

from (\ref{f_3}) we obtain the relation
\begin{eqnarray}
u_j=\frac{1}{\left(e_j,e_{j-1}^\bot\right)\sqrt{l_{j-1}l_j}}v_{j-1}\label{f_4}\\
+\frac{1}{l_j}\left(\frac{\left(e_{j+1},e_{j}\right)}{\left(e_{j+1},e_{j}^\bot\right)}-\frac{\left(e_{j},e_{j-1}\right)}{\left(e_{j}^\bot,e_{j-1}\right)}\right)v_j
-\frac{1}{\left(e_j^\bot,e_{j+1}\right)\sqrt{l_jl_{j+1}}}v_{j+1}.\notag
\end{eqnarray}
On introducing the notations
\begin{eqnarray*}
\rho_j=\frac{-1}{\left(e_{j+1},e_{j}^\bot\right)\sqrt{l_{j}l_{j+1}}},\quad
j\geqslant 1,\\
q_j=\frac{1}{l_j}\left(\frac{\left(e_{j},e_{j+1}\right)}{\left(e_{j}^\bot,e_{j+1}\right)}-\frac{\left(e_{j},e_{j-1}\right)}{\left(e_{j}^\bot,e_{j-1}\right)}\right),\quad
j\geqslant 2,
\end{eqnarray*}
we can rewrite (\ref{f_4}) in a form:
\begin{equation*}
u_j=\rho_{j-1}v_{j-1}+q_jv_j+\rho_jv_{j+1}, \quad j\geqslant 2,
\end{equation*}
and $q_1$ is found from the condition at zero. So finally we
obtain the following result: if $f$ and $g$ having representation
(\ref{F_form}) are connected by (\ref{EqnCS}), then $u$ and $v$
defined by (\ref{Subst}) satisfy
\begin{equation*}
Av=u,\quad A=\begin{pmatrix} q_1& \rho_1& 0& 0& 0\\
\rho_1& q_2 &\rho_2 &0&0\\
0&\rho_2&q_3&\rho_3&0\\
0 &0 & \cdot&\cdot&\cdot
\end{pmatrix}
\end{equation*}
We can introduce the dependence on (continuous) time $t$: let
$f(x,t),$ $g(x,t)$ have form:
\begin{eqnarray*}
f(x,t)=f_j(t)e_j(x)+\xi(x,t)e_j^\bot(x),\quad
x\in \Delta_j, \\
g(x,t)=g_j(t)e_j(x)+\eta(x,t)e_j^\bot(x),\quad x\in \Delta_j,
\end{eqnarray*}
then if $g(x,t)=if_t(x,t)$, then $f$ solves
\begin{equation*}
Jf_x=iHf_t.
\end{equation*}
On the other hand (\ref{Subst}) implies the relationship
$u_j(t)=i{v_j}_t(t)$, which yields that $v$ solves $iv_t -Av=0.$
Adding initial and boundary conditions gives well-posed IBVP for
dynamical system with continuous time governed by Jacobi matrix:
\begin{equation}
\label{Jac_SchrDyn} \left\{
\begin{array}l
iv_t -Av=0,\quad x\geqslant 0,\,t\geqslant 0,\\
v_n(0)=0,\quad n\geqslant 1,\\
v_1(t)=h(t),\quad t\geqslant 0.
\end{array}
\right.
\end{equation}
The response operator $R^T_J:L_2(0,T)\mapsto L_2(0,T)$ with the
domain $D=C_0^\infty(0,T)$ for this system is introduced by the
rule $\left(R^T_Jh\right)(t):=v_2(t)$. On the other hand, IBVP
(\ref{Jac_SchrDyn}) is equivalent to (we assume that
$e_1=(1,0)^T$):
\begin{equation}
\label{Jac_Can_syst} \left\{
\begin{array}l
iHf_t-Jf_x=0,\quad x\geqslant 0,\,t\geqslant 0,\\
f(x,0)=0,\quad x\geqslant 0,\\
f^1(0,t)=j(t):=\frac{h(t)}{\sqrt{l_1}},\quad t\geqslant 0.
\end{array}
\right.
\end{equation}
For the system (\ref{Jac_Can_syst}) the response operator
$\widetilde R^T_J:L_2(0,T)\mapsto L_2(0,T)$ with the domain
$D=C_0^\infty(0,T)$ is introduced by the rule
$\left(R^T_Jh\right)(t):=f^2(0,t)$. Note that by (\ref{F_form}),
$f^2(0,t)=\xi_1(0,t)$. From (\ref{Xi_expr}), the relationship
$g(x,t)=if_t(x,t)$ and (\ref{f_1}), we have that
\begin{eqnarray*}
\left(R^T_{J}h\right)(t)=s_1(t)+g_1(t)l_1=\frac{f_{2}(t)}{\left(e_2,e_{1}^\bot\right)}-f_{1}(t)+if^1(0,t)l_1\\=
\frac{f_{2}(t)}{\left(e_2,e_{1}^\bot\right)}-\frac{h(t)}{\sqrt{l_1}}+ih(t)\sqrt{l_1}=
-\rho_1v_{2}(t)\sqrt{l_1}-h(t)\left(\frac{1}{\sqrt{l_1}}-i\sqrt{l_1}\right).
\end{eqnarray*}
So IP for (\ref{Jac_SchrDyn}) and (\ref{Jac_Can_syst}) from
corresponding response operators are equivalent. We note that that
we can introduce the different type of  continuous dynamics for
Jacobi matrices (for example the dynamics of the type
$\frac{d}{dt^2}$ was considered in \cite{T}).

We can also introduce the dependence on the discrete time $t\in
\mathbb{N}$ by letting $f_t(x),$ $g_t(x)$ have form:
\begin{eqnarray*}
f_t(x)=f_{j,t}e_j(x)+\xi_t(x)e_j^\bot(x),\quad
x\in \Delta_j,\,t\in \mathbb{N}, \\
g_t(x)=g_{j,t}e_j(x)+\eta_t(x)e_j^\bot(x),\quad x\in
\Delta_j,\,t\in \mathbb{N}.
\end{eqnarray*}
If $f,g$ are related by
$g_t(x)=f_t(x)+f_{t-1}(x)=:\partial_tf(x)$, then counting
(\ref{EqnCS}), $f$ solves
\begin{equation*}
Jf_x=H\partial_t f.
\end{equation*}
The equality (\ref{Subst}) implies $u_j=\partial_tv_j,$ which
yields that $v$ satisfies $\partial_t v_{\cdot,t}-Av_{\cdot,t}=0$.
Adding initial and boundary conditions gives the following IBVP:
\begin{equation}
\label{Jac_SchrDyn_discr} \left\{
\begin{array}l
\partial_t v_{\cdot,t}-Av_{\cdot,t}=0,\quad t\in \mathbb{N}\\
v_{n,1}=v_{n,0}(0)=0,\quad n\geqslant 1,\\
v_{1,t}=h_t,\quad t\in \mathbb{N}.
\end{array}
\right.
\end{equation}
where $h_t\in l_2$ is referred to as a \emph{boundary control}.
The response operator $R^T_{J,d}$ with the domain $D=\mathrm{R}^T$
for this system is introduced by $R^T_{J,d}: \mathbb{R}^T\mapsto
\mathbb{R}^T$, $\left(R^T_{J,d}h\right)_t=v_{2,t},$ $t=1\ldots,T$.
The forward and inverse problem was studied in \cite{MM,MM1}. The
IBVP (\ref{Jac_SchrDyn_discr}) is equivalent to which is
equivalent to the following IBVP for a canonical system:
\begin{equation}
\label{Jac_Can_syst_dscr} \left\{
\begin{array}l
H\partial_t f-Jf_x=0,\quad x\geqslant 0, \,t\in \mathbb{N},\\
f_0(x)=0,\quad x\geqslant 0,\\
f^1_t(0)=j_t:=\frac{h_t}{\sqrt{l_1}},\quad t\in \mathbb{N}.
\end{array}
\right.
\end{equation}
For the system (\ref{Jac_Can_syst_dscr}) the response operator
$\widetilde R^T_{J,d}:l_2\mapsto l_2$ is introduced by the rule
$\left(R^T_Jj\right)(t):=f^2_t(0)$. By (\ref{F_form}),
$f^2_t(0)={\xi_1}_t(0)$, from (\ref{Xi_expr}), the relationship
$g_t(x)=\partial_tf(x)$ and (\ref{f_1}), we have that
\begin{eqnarray*}
\left(R^T_{J,d}h\right)_t={s_1}_t+{g_1}_tl_1=\frac{{f_{2}}_t}{\left(e_2,e_{1}^\bot\right)}-{f_{1}}_t+if^1_t(0)l_1\\=
\frac{{f_{2}}_t}{\left(e_2,e_{1}^\bot\right)}-\frac{h_t}{\sqrt{l_1}}+ih_t\sqrt{l_1}=
-\rho_1v_{2,t}\sqrt{l_1}-h_t\left(\frac{1}{\sqrt{l_1}}-i\sqrt{l_1}\right).
\end{eqnarray*}
So IP for (\ref{Jac_SchrDyn_discr}) and (\ref{Jac_Can_syst_dscr})
from corresponding response operators are equivalent.

We see that different dynamic systems after transformations come
to dynamical canonical systems with different dynamics
($i\frac{d}{dt}$, $\frac{d}{dt^2}$, and even discrete one
$\partial_t$).

We will investigate the dynamics given by $i\frac{d}{dt}$, the
canonical system with this dynamics possess property of finite
speed of wave propagation.

\section{Canonical systems with smooth strictly positive Hamiltonian.}

We consider the IBVP for a canonical system. Assuming that the
Hamiltonian satisfies conditions: $H=H^*\in
C^2(0,T;\mathbb{R}^{2\times 2})$, $H\geqslant \delta>0$,
$\operatorname{tr}{H}=1$,  we set
$Y^f=\begin{pmatrix}y^1\\y^2\end{pmatrix}$ to be a solution to
\begin{equation}
\label{CanSyst0} \left\{
\begin{array}l
iH\frac{d}{dt}Y-J\frac{d}{dx}Y=0,\quad x\geqslant 0,\,t\geqslant 0,\\
Y(x,0)=0,\quad x\geqslant 0,\\
y^1(0,t)=f(t),\quad t\geqslant 0.
\end{array}
\right.
\end{equation}
Where the \emph{boundary control} $f\in
\mathcal{F}^T:=L_2(0,T;\mathbb{C}).$ The \emph{response operator}
$R^T:\mathcal{F}^T\mapsto \mathcal{F}^T$ is introduced as
$\left(R^Tf\right)(t):=y_2^f(0,t)$. The inverse problem we will be
dealing with consists in a recovering $H(x),$ on an interval
$(0,l)$ for some $l>0$ from given $R^{2T}$.

\subsection{One-velocity wave system.}

We rewrite (\ref{CanSyst0}): differentiate the first line in
(\ref{CanSyst0}) w.r.t. $t$ and use equation to get:
\begin{equation*}
HY_{tt}+JH^{-1}JY_{xx}+JH_x^{-1}JY_x=0,
\end{equation*}
which is equivalent to the equation
\begin{equation*}
HY_{tt}-\frac{1}{\det{H}}HY_{xx}+JH_x^{-1}JY_x=0.
\end{equation*}
Counting the initial and boundary condition, we obtain that $Y$
satisfies the following IBVP for one-velocity system:
\begin{equation}
\label{WaveEqn0} \left\{
\begin{array}l
\det{H}Y_{tt}-Y_{xx}+\det{H}H^{-1}JH_x^{-1}JY_x=0,\quad x\geqslant 0,\, t\geqslant 0,\\
Y(x,0)=Y_t(x,0)=0,\quad x\geqslant 0,\\
\begin{pmatrix}y^1(0,t) \\
y^2(0,t)\end{pmatrix}=G(t):=\begin{pmatrix}f(t)\\(Rf)(t)\end{pmatrix},\quad
t\geqslant 0.
\end{array}
\right.
\end{equation}
Here the velocity is given by $c(x)=\frac{1}{\sqrt{\det{H(x)}}}$.
The response operator $R^T_w:L_2(0,T;\mathbb{C})\mapsto
L_2(0,T;\mathbb{C})$ with the domain
$\mathcal{D}=C_0^\infty(0,T,\mathbb{C})$ for (\ref{WaveEqn0}) is
introduced as , $\left(R^T_wG\right)(t):=Y^G_x(0,t)$. The
\emph{eikonal} function is introduced by
$\tau(x):=\int_0^x\sqrt{\det H(s)}\,ds,$ and $\Omega^l=\{x>0\,|\,
\tau(x)<l\}$. Then the natural setup of IP is to recover
$H(x)|_{\Omega^T}$ from $R^{2T}_w$.

We see that the IP for the system (\ref{WaveEqn0}), is equivalent
to IP for (\ref{CanSyst0}). But there is one important disadvantage -- in studying IP for (\ref{WaveEqn0})
which comes from (\ref{CanSyst0}), we need to use the specific set of controls of the type $\begin{pmatrix} f \\
Rf\end{pmatrix}$, which makes application of the BC method
problematic. Instead, we will reduce (\ref{CanSyst0}) to
Dirac-type system, and follow the scheme offered in \cite{BM_Dir}.

\subsection{Dirac-type dynamical system.}

We introduce the following transformation: let
\begin{equation*}
U=\begin{pmatrix} \cos{\phi(x)} & \sin{\phi(x)}\\
-\sin{\phi(x)} & \cos{\phi(x)}\end{pmatrix}
\end{equation*}
be a unitary matrix such that $U^*HU=D:=\begin{pmatrix} d_1(x) & 0\\
0 & d_2(x)\end{pmatrix}$, where $d_1,d_2\geqslant\delta>0$,
$d_1+d_2=1.$ If $Y=U\widetilde Y$, then $\widetilde Y$ satisfies
the following IBVP for Dirac-type dynamical system:
\begin{equation}
\label{Can_syst_Dir} \left\{
\begin{array}l
iD\frac{d}{dt}\widetilde Y+J\frac{d}{dx}\widetilde
Y-\phi'(x)\widetilde Y=0,\quad x\geqslant 0,\, t\geqslant 0,\\
\widetilde Y(x,0)=0,\quad x\geqslant 0,\\
\widetilde
y^1(0,t)=g(t):=\cos{\phi(0)}f(t)+\sin{\phi(0)}(Rf)(t),\quad
t\geqslant 0.
\end{array}
\right.
\end{equation}
The response operator $R^T_{CD}:L_2(0,T)\mapsto L_2(0,T)$ is
introduced by $\left(R^T_{CD}g\right)(t):=\widetilde y^2(0,t)$. We
can see that $\widetilde
y^2(0,t)=-\sin{\phi(0)}f(t)+\cos{\phi(0)}(Rf)(t),$ so IP for
(\ref{CanSyst0}) and for (\ref{Can_syst_Dir}) are equivalent.

Thus our first goal will be to study the dynamic IP for the
following Dirac-type system:
\begin{equation}
\label{DiracSyst}
\left\{
\begin{array}l
iD\frac{d}{dt}V+J\frac{d}{dx}
V+\psi(x)V=0,\quad x\geqslant0, \,t\geqslant 0,\\
V(x,0)=0,\quad x\geqslant0,\\
v^1(0,t)=f(t), \quad t\geqslant 0,
\end{array}
\right.
\end{equation}
where $D$ as above is a diagonal matrix with twice differentiable
entries and unit trace, $\psi\in C^2(\mathbb{R}_+)$. The function
$f\in \widetilde{\mathcal{F}}^T:=L_2(0,T;\mathbb{C})$ is a
\emph{boundary control}. The response
$R^T_D:\widetilde{\mathcal{F}}^T\mapsto \widetilde{\mathcal{F}}^T$
is introduced by $\left(R^T_Df\right)(t):=v^2(0,t)$. The IP
consists in recovering $D|_{\Omega^T},\,\psi|_{\Omega^T}$ from
$R^{2T}$. We outline the scheme offered in \cite{BM_Dir,MM2}:
\begin{proposition}
\label{Repr} The solution to (\ref{DiracSyst}) admits the
following representation:
\begin{equation*}
V(x,t)=A(x)f(t-\tau(x))+\int_0^{x(t)}w(x,s)f(t-\tau(s))\,ds,
\end{equation*}
where $\tau(s)=\int_0^s\sqrt{d_1(\alpha)d_2(\alpha)}\,d\alpha$ is
eikonal, $x(t)$ is a function inverse to $\tau(x)$, the kernel
$w=\begin{pmatrix}w^1\\w^2\end{pmatrix}$ is twice differentiable
in $\left\{(x,s)\,|\, 0\leqslant\tau(x)\leqslant s\leqslant
T\right\}$, $A=\begin{pmatrix}a^1\\a^2\end{pmatrix}$, where
$a^{1,2}$ are solutions to the following system
\begin{eqnarray*}
i\sqrt{d_1}a^1_x=\sqrt{d_2}a^2_x,\\
\sqrt{d_2}\left(\psi a^1+a^2_x\right)=i\sqrt{d_1}\left(\psi
a^2-a^1_x\right).
\end{eqnarray*}

\end{proposition}
We introduce the \emph{outer space}, the space o states of
(\ref{DiracSyst}): $\mathcal{H}^T:=L_2(0,\tau(T);\mathbb{C})$ and
a \emph{control operator} $\widetilde
W^T:\widetilde{\mathcal{F}}^T\mapsto \mathcal{H}^T$ acting by the
rule
\begin{equation*}
\left(\widetilde W^Tf\right)(x):=V^f(x,T).
\end{equation*}
The Proposition \ref{Repr} implies that $\widetilde W^T$ is not an
isomorphism, and the system (\ref{DiracSyst}) is not boundary
controllable. To restore the controllability, we introduce the
auxiliary system:
\begin{equation}
\label{DiracSyst1} \left\{
\begin{array}l
iD\frac{d}{dt}U-J\frac{d}{dx}
U-\psi(x)U=0,\quad x\geqslant0, \,t\geqslant 0,\\
U(x,0)=0,\quad x\geqslant0,\\
u^1(0,t)=g(t), \quad t\geqslant 0,
\end{array}
\right.
\end{equation}
and note that solutions to (\ref{DiracSyst}) and
(\ref{DiracSyst1}) are connected by the formula
$V^f=\overline{U^{\overline f}}$. The extended outer space is
defined by $\mathcal{F}^T:=L_2(0,T;\mathbb{C}^2)$, and the
\emph{extended control operator} $W^T: \mathcal{F}^T\mapsto
\mathcal{H}^T$ is introduced by
\begin{equation*}
W^T\begin{pmatrix} f\\g\end{pmatrix}:=V^f(x,T)+U^g(x,T).
\end{equation*}
\begin{proposition}
\label{PropContr} The extended control operator is an isomorphism
between $\mathcal{F}^T$ and $\mathcal{H}^T$.
\end{proposition}
The set $\mathcal{U}^T:=W^T\mathcal{F}^T$ is called extended
reachable set. The Proposition \ref{PropContr} says that
$\mathcal{U}^T=\mathcal{H}^T$.

We consider the operator of the Dirac-type system on a half-line:
let $\mathbf{D}:=D^{-1}J\frac{d}{dx}+D^{-1}\psi$ on
$L_2(\mathbb{R}_{+},\mathbb{C}^2)\ni\Phi=\begin{pmatrix}\Phi_1\\\Phi_2\end{pmatrix}$
with a Dirichlet condition $\Phi_1(0)=0$. Denote by
$\theta(x,z)=\begin{pmatrix}\theta_1\\ \theta_2\end{pmatrix}$ a
solution to the following Cauchy problem for $z\in\mathbb{C}$:
\begin{equation}
\label{spec_sol} \left\{\begin{array}l
J\theta_x+V\theta=zD\theta, \quad x>0, \\
\theta_1(0,z)=0, \quad \theta_2(0,z)=1.
\end{array}\right.
\end{equation}
Let $d\rho$ be a spectral measure of $\mathbf{D}$, and $F:
L_2(\mathbb{R}_+;\mathbb{C}^2)\mapsto L_{2,\,\rho}(\mathbb{R}_+)$
be the corresponding Fourier transform:
\begin{eqnarray*}
\left(F\begin{pmatrix}f_1\\
f_2\end{pmatrix}\right)(\lambda)=F(\lambda)=\int_0^\infty\left(
f_1(x)\theta_1(x,\lambda)+f_2(x)\theta_2(x,\lambda)\right)\,dx,\\
f_1(x)=\int_{-\infty}^\infty
F(\lambda)\theta_1(x,\lambda)\,d\rho(\lambda),\,\,
f_2(x)=\int_{-\infty}^\infty F(\lambda)\theta_2(x,\lambda)\,d\rho(\lambda),\\
\int_0^\infty\left(
f_1^2(x)+f_2^2(x)\right)\,dx=\int_{-\infty}^\infty
F^2(\lambda)\,d\rho(\lambda).
\end{eqnarray*}
We introduce the \emph{extending connecting operator}
$C^T:\mathcal{F}^T\mapsto \mathcal{F}^T$ by the quadratic form:
\begin{equation}
\label{CT}
\left(C^T\begin{pmatrix} f_1 \\
g_1\end{pmatrix},\begin{pmatrix} f_2 \\
g_2\end{pmatrix}\right)_{\mathcal{F}^T}=\left(W^T\begin{pmatrix} f_1 \\
g_1\end{pmatrix},W^T\begin{pmatrix} f_2 \\
g_2\end{pmatrix}\right)_{\mathcal{H}^T}, \,\,
C^T=\left(W^T\right)^*W^T.
\end{equation}
The important fact in the BC method is that
\begin{proposition}
The extending connecting operator is a positive isomorphism in
$\mathcal{F}^T$, it admits the representation in terms of dynamic
inverse data $R^{2T}$, and spectral inverse data $d\rho(\lambda)$.
\end{proposition}
We introduce the linear manifold of Fourier images of extended
states (Fourier image of extended reachable set) at time $t=T$:
\begin{equation*}
B^T_D:=\left\{ K(\lambda)\,\bigl|\, K(\lambda)=\left(FW^T\begin{pmatrix}k_1 \\
k_2\end{pmatrix}\right)(\lambda),\,
\begin{pmatrix}k_1\\k_2\end{pmatrix}\in
\mathcal{F}^T\right\}=F\mathcal{U}^T,.
\end{equation*}
Equipped with the scalar product, generated by $C^T$:
\begin{equation*}
[F,G]_{B^T_D}:=\left(C^T\begin{pmatrix}f_1 \\
f_2\end{pmatrix},\begin{pmatrix}g_1 \\
g_2\end{pmatrix}\right)_{\mathcal{F}^T},\quad F,G\in B^T_D,\\
\end{equation*}
this linear space becomes a Hilbert space of analytic functions.
It is also possible to define a reproducing kernel in this space
(it is given in terms of a solution to a Krein equation), which
makes $B^T_D$ a de Branges space. Solution of dynamic and spectral
IPs for (\ref{DiracSyst}) and construction of corresponding de
Branges space will be the subject of forthcoming publications.

\subsection{Dynamic approach to de Branges spaces.}

Basing on the arguments from the previous subsection, we can
formulate the hypothesis about de Branges space for canonical
system (\ref{CanSyst0}) with general Hamiltonian. First we
introduce the auxiliary system
\begin{equation}
\label{CanSyst1} \left\{
\begin{array}l
iH\frac{d}{dt}Z+J\frac{d}{dx}Z=0,\quad x\geqslant 0,\,t\geqslant 0,\\
Z(x,0)=0,\quad x\geqslant 0,\\
z^1(0,t)=g(t),\quad t\geqslant 0.\\
\end{array}
\right.
\end{equation}
The \emph{extending control operator} $W^T: \mathcal{F}^T\mapsto
\mathcal{F}^T$ acting in extended control space
$\mathcal{F}^T:=L_2(0,T;\mathbb{C}^2)$ is defined by
$W^T\begin{pmatrix} f\\g\end{pmatrix}:=Y^f(x,T)+Z^g(x,T).$ The
\emph{extending connecting operator} $C^T$ is given by analog to
(\ref{CT}). Then the de Branges space corresponding to
(\ref{CanSyst0}) is a Fourier image of extended reachable set,
equipped with a scalar product, generated by $C^T$.

We note that the construction of de Branges space by dynamic
methods for general Hamiltonian in fact is equivalent to solving
the dynamic IP for system (\ref{CanSyst0}) with general $H$. We
note that the in studying the IP in this case, one inevitably face
with two obstacles: the smoothness of $H$, and changing the rank
of $H$, which reflects in the lack of the boundary controllability
of the dynamical system. The authors suggest that studying the
inverse dynamic problem for a Krein string \cite{Kr2,DMcK} will be
instructive and can help to overcome difficulties connected with
general Hamiltonian.

\noindent{\bf Acknowledgments}

The research of Victor Mikhaylov was supported in part by RFBR
17-01-00529-A. Alexandr Mikhaylov was supported by RFBR
17-01-00099-A; A. S. Mikhaylov and V. S. Mikhaylov were partly
supported by RFBR 18-01-00269-A and VW Foundation program
"Modeling, Analysis, and Approximation Theory toward application
in tomography and inverse problems." The authors are deeply
indebted to Prof. R. V. Romanov and Prof. M. I. Belishev for
valuable discussions.


\begin{thebibliography}{99}

\bibitem{ALP1}
\textsc{S. A. Avdonin, S. Lenhart and V. Protopopescu.}
Determining the potential in the Schr\"{o}dinger equation from the
Dirichlet to Neumann map by the boundary control method \textit{J.
Inverse Ill-Posed Probl,} 13, 2005, 317--330.


\bibitem{AM}
\textsc{S. A. Avdonin, V. S. Mikhaylov}. The boundary control
approach to inverse spectral theory. \textit{Inverse Problems,}
26, no. 4, 2010, 045009, 19 pp.

\bibitem{B87} \textsc{M. I. Belishev.}  An approach to multidimensional inverse problems for the wave equation
\textit{Dokl. Akad. Nauk SSSR,} 297, no. 3, 1987, 524--527 (in
Russian)
\newline \textsc{M. I. Belishev.} \textit{Soviet Math. Dokl.} 36,
no. 3, 1988, 481--484 (Engl. Transl.)

\bibitem{B97} \textsc{M. I. Belishev.} Boundary control in reconstruction of manifolds
and metrics (the BC method). \textit{Inverse Problems} 13, no. 5,
1997, R1--R45.

\bibitem{B08}
\textsc{M. I. Belishev.} Boundary control and inverse problems: a
one-dimensional version of the boundary control method. \textit{
Zap. Nauchn. Sem. S.-Peterburg. Otdel. Mat. Inst. Steklov.
(POMI),} 354, 2008, 19--80 (in Russian); English translation:
\textit{J. Math. Sci. (N. Y.),} 155, no 3, 2008, 343-v378.

\bibitem{BIP07}
\textsc{M.I.Belishev}. Recent progress in the boundary control
method. \textit{Inverse Problems}, 23, no 5, 2007, R1--R67.

\bibitem{B17} \textsc{M.I.Belishev.} Boundary control and tomography of Riemannian
manifolds (the BC-method). \textit{Uspekhi Matem. Nauk,} 72, no 4,
2017, 3-66 (in Russian).

\bibitem{BM_Sh}
\textsc{M.I.Belishev, V.S.Mikhaylov}. Unified approach to
classical equations of inverse problem theory. \textit{Journal of
Inverse and Ill-Posed Problems}, 20, no 4, 2012, 461--488.

\bibitem{BM_Dir}
\textsc{M. I. Belishev, V. S. Mikhaylov}. Inverse problem for
one-dimensional dynamical Dirac system (BC-method).
\textit{Inverse Problems}, 26, no. 4, 2010, 045009, 19 pp.

\bibitem{DBr}
\textsc{Louis de Branges}. Hilbert space of entire functions.
\textit{Prentice-Hall, NJ.} 1968.

\bibitem{DMcK}
\textsc{H. Dym, H. P. McKean}. Gaussian processes, function
theory, and the inverse spectral problem. \textit{Academic Press,
New York etc.} 1976.

\bibitem{Kr2} \textsc{M. G. Krein.} On the one method of effective solving the inverse boundary
value problem \textit{Dokl. Akad. Nauk. SSSR}, 94, no. 6, 1954,
987--990.


\bibitem{MM}
\textsc{A. S. Mikhaylov, V. S Mikhaylov}. Dynamical inverse
problem for the discrete Schr\"odinger operator.
\textit{Nanosystems: Physics, Chemistry, Mathematics,} 7, no. 5,
2016, 842-854.

\bibitem{MM3}
\textsc{A. S. Mikhaylov, V. S Mikhaylov}. Relationship between
different types of inverse data for the one-dimensional
Schr\"odinger operator on a half-line. \textit{Zapiski Nauchnykh
Seminarov POMI}, 451, 2016, 134-155.

\bibitem{MM1}
\textsc{A. S. Mikhaylov, V. S Mikhaylov}. Dynamic inverse problem
for the Jacobi matrices.
\textit{https://arxiv.org/abs/1704.02481}.

\bibitem{MM2}
\textsc{A. S. Mikhaylov, V. S Mikhaylov}. The Boundary Control
method and de Branges spaces. Schr\"odinger equation, Dirac system
and Discrete Schr\"odinger operator. \textit{Journal of
Mathematical Analysis and Applications}, doi:
10.1016/j.jmaa.2017.12.013, 2017.

\bibitem{R2} \textsc{C. Remling }. Schr\"odinger operators and de Branges spaces. \textit{J.
Funct. Anal.} 196, no. 2, 2002, 323--394.

\bibitem{RR} \textsc{R. V. Romanov.} Canonical systems and de Branges spaces
\textit{http://arxiv.org/abs/1408.6022}.

\bibitem{T} \textsc{ G. Teschl } Jacobi operators and completely integrable
nonlinear lattices. Mathematical Surveys and Monographs, 72.
\textit{American Mathematical Society, Providence, RI,} 2000.
\end{thebibliography}
\end{document}